\documentclass[10pt]{article}
\usepackage{setspace}

\usepackage{amsmath} 
\usepackage{amssymb}
\usepackage{latexsym}
\usepackage{xcolor}
\usepackage{fancyhdr}
\allowdisplaybreaks 

 \usepackage{comment}

\def\XXint#1#2#3{{\setbox0=\hbox{$#1{#2#3}{\int}$} 
\vcenter{\hbox{$#2#3$}}\kern-.5\wd0}}   

 \numberwithin{equation}{section}
\newtheorem{thm}[equation]{Theorem}
\newtheorem{prop}[equation]{Proposition}
\newtheorem{defn}[equation]{Definition}

\newtheorem{lem}[equation]{Lemma}

\title{A necessary condition on a singular kernel for the continuity of an integral operator in H\"{o}lder spaces}

 \author{  
Massimo Lanza de Cristoforis
\\
Dipartimento di Matematica `Tullio Levi-Civita', 
\\
Universit\`a degli Studi di Padova, 
\\
Via Trieste 63, Padova 35121, 
Italy. 
\\
E-mail: mldc@math.unipd.it   }

\date{\ }
 
\begin{document}
 
 \maketitle

\noindent
{\bf Abstract:} We prove that a condition of boundedness of the maximal function of a singular integral operator, that is known to be sufficient for the continuity of a corresponding integral operator in H\"{o}lder spaces, is actually also necessary  in case the action of the integral operator does not decrease the regularity of a function. We do so in the frame of  metric measured spaces with a measure satisfying  certain growth conditions that include  nondoubling measures. Then we present an application to the case of an integral operator defined on a compact differentiable manifold. 

 \vspace{\baselineskip}





\section{Introduction}
This paper concerns the analysis of a singular integral operator   in H\"{o}lder spaces on subsets of measured metric spaces with a measure that may be non-doubling, in the spirit of a research project that has been initiated by Garc\'{\i}a-Cuerva and  Gatto, \cite{GaGa04}, \cite{GaGa05}, Gatto \cite{Gat06}, \cite{Gat09} and that is here being developed with the aim of specific applications to potential theory. Let 
$
(M,d)$ be a metric space and let 
$X$, $Y$ be subsets of $M$. 
\begin{eqnarray} \nonumber
&&\text{Let}\ {\mathcal{N}} \ \text{ be a $\sigma$-algebra of parts of}\  Y  \,, {\mathcal{B}}_Y\subseteq  {\mathcal{N}}\,.
\\   \label{eq:nu}
&&\text{Let}\ \nu\  \text{be   measure on}\  {\mathcal{N}} \,.
\\  \nonumber
&&\text{Let}\ \nu(B(x,r)\cap Y)<+\infty\qquad\forall (x,r)\in X\times ]0,+\infty[\,,
\end{eqnarray}
where  ${\mathcal{B}}_Y$ denotes the $\sigma$-algebra of the Borel subsets of $Y$ and 
\begin{equation}\label{eq:balls}
B(\xi,r)\equiv \left\{\eta\in M:\, d(\xi,\eta)<r\right\}\,,\quad
B(\xi,r]\equiv \left\{\eta\in M:\, d(\xi,\eta)\leq r\right\}
\,,
\end{equation}
for all $(\xi,r)\in M\times ]0,+\infty[$. We assume that  $\upsilon_Y\in ]0,+\infty[$ and we consider two types of assumptions on $\nu$. 
The first assumption is  that
$Y$ is   upper $\upsilon_Y$-Ahlfors regular with respect to $X$, \textit{i.e.}, that
\begin{eqnarray} \nonumber
&&\text{there\ exist}\ r_{X,Y,\upsilon_Y}\in]0,+\infty]\,,\ c_{X,Y,\upsilon_Y}\in]0,+\infty[\ \text{such\ that}
\\ \nonumber
&&\nu( B(x,r)\cap Y )\leq c_{X,Y,\upsilon_Y} r^{\upsilon_Y}  
\\  \label{defn:uareg1}
&&\text{for\ all}\ x\in X\ \text{and}\  r \in]0,r_{X,Y,\upsilon_Y}[
 \,. 
\end{eqnarray}
In case $X=Y$, we just say that $Y$ is   upper $\upsilon_Y$-Ahlfors regular  and this is the assumption that has been considered by 
Garc\'{\i}a-Cuerva and Gatto \cite{GaGa04}, \cite{GaGa05}, Gatto \cite{Gat06}, \cite{Gat09}   in case $X=Y=M$. See also  Edmunds,  Kokilashvili and   Meskhi~\cite[Chap.~6]{EdKoMe02} in the frame of Lebsgue spaces.
  
 Then we consider a stronger version of the upper Ahlfors regularity as in \cite{La22a}. Namely, 
  we assume that $Y$ is strongly upper $\upsilon_Y$-Ahlfors regular with respect to $X$, \textit{i.e.}, that
\begin{eqnarray} \nonumber
&&\text{there\ exist}\ r_{X,Y,\upsilon_Y}\in]0,+\infty]\,,\ c_{X,Y,\upsilon_Y}\in]0,+\infty[\ \text{such\ that}
\\ \nonumber
&&\nu( (B(x,r_2)\setminus B(x,r_1))\cap Y )\leq c_{X,Y,\upsilon_Y}(r_2^{\upsilon_Y}-r_1^{\upsilon_Y})
\\ \label{defn:suareg1}
&&\text{for\ all}\ x\in X\ \text{and}\ r_1,r_2\in[0,r_{X,Y,\upsilon_Y}[
\ \text{with}\ r_1<r_2\,,
\end{eqnarray}
where we understand that $B(x,0)\equiv\emptyset$  (in case $X=Y$, we just say that $Y$ is strongly upper $\upsilon_Y$-Ahlfors regular). So, for example,	
 if $Y$ is the boundary of a bounded open Lipschitz  subset of $M={\mathbb{R}}^n$, then $Y$ is upper $(n-1)$-Ahlfors regular with respect to ${\mathbb{R}}^n$ 
and if $Y$ is the boundary of an open bounded subset of $M={\mathbb{R}}^n$ of class $C^1$, then $Y$ is strongly upper $(n-1)$-Ahlfors regular  with respect to $Y$. Here and throughout the paper,  
\[
n\in {\mathbb{N}}, n\geq 2\,.
\]
We plan to consider  off-diagonal kernels $K$ from $(X\times Y)\setminus D_{X\times Y}$ to ${\mathbb{C}}$, where
\[
D_{X\times Y}\equiv  \left\{
(x,y)\in X\times Y:\,x=y
\right\}
\] denotes the diagonal set of $X\times Y$  
 and we introduce the following class of   `potential type' kernels  (see also   \cite{DoLa17}).	 
\begin{defn}\label{defn:ksss} Let 
$
(M,d)$ be a metric space.
 Let $X$, $Y\subseteq M$. Let $s_1$, $s_2$, $s_3\in {\mathbb{R}}$. We denote by the symbol ${\mathcal{K}}_{s_1, s_2, s_3} (X\times Y)$ the set of continuous functions $K$ from $(X\times Y)\setminus D_{X\times Y}$ to ${\mathbb{C}}$ such that
 \begin{eqnarray*}
\lefteqn{
\|K\|_{  {\mathcal{K}}_{ s_1, s_2, s_3  }(X\times Y)  }
\equiv
\sup\biggl\{\biggr.
d(x,y)^{ s_{1} }\vert K(x,y)\vert :\,(x,y)\in X\times Y, x\neq y
\biggl.\biggr\}
}
\\ \nonumber
&&\qquad\qquad\qquad
+\sup\biggl\{\biggr.
\frac{d(x',y)^{s_{2}}}{d(x',x'')^{s_{3}}}
\vert   K(x',y)- K(x'',y)  \vert :\,
\\ \nonumber
&&\qquad\qquad\qquad 
x',x''\in X, x'\neq x'', y\in Y\setminus B(x',2d(x',x''))
\biggl.\biggr\}<+\infty\,.
\end{eqnarray*}
\end{defn}
For  $s_2=s_1+s_3$ one has the so-called class of standard kernels  that is the case in which Garc\'{\i}a-Cuerva and Gatto \cite{GaGa04}, \cite{GaGa05}, Gatto \cite{Gat09} have proved  $T1$ Theorems for the integral operators with kernel $K$ in case of weakly singular, singular and hyper-singular integral operators with $X=Y$.  We plan to analyze the integral operator
\begin{equation}\label{prop:qz6}
Q[Z,g,1](x)\equiv \int_YZ(x,y)(g(x)-g(y))\,d\nu(y)\qquad\forall x\in X\,.
\end{equation}
where  $Z$ belongs to a class ${\mathcal{K}}_{ s_1, s_2, s_3  } (X\times Y)$ as in Definition \ref{defn:ksss}
and $g $ is a ${\mathbb{C}}$-valued function  in  $X\cup Y$. We exploit the operator in (\ref{prop:qz6}) in the analysis of the double layer potential and we note that operators as in (\ref{prop:qz6})  appear in the applications (cf. \textit{e.g.},   \cite[\S 8]{DoLa17}). 
As in (cf.~\cite[Prop.~6.3 (ii)]{La22a}), we  can estimate the H\"{o}lder quotient of $Q[Z,g,1]$ by introducing  a more restrictive class of kernels. Namely the following.
\begin{defn}Let 
$
(M,d)$ be a metric space. Let $X$, $Y\subseteq M$. Let $\nu$ be as in (\ref{eq:nu}).   Let $s_1$, $s_2$, $s_3\in {\mathbb{R}}$. We set
\begin{eqnarray*}\nonumber
\lefteqn{
 {\mathcal{K}}_{ s_1, s_2, s_3  }^\sharp(X\times Y)
  \equiv 
  \biggl\{\biggr.
K\in  {\mathcal{K}}_{ s_1, s_2, s_3  }(X\times Y):\,
}
\\ \nonumber
&&\ \ 
K(x,\cdot)\ \text{is}\ \nu-\text{integrable\ in}\ Y\setminus B(x,r)
\ \text{for\ all}\ (x,r)\in X\times]0,+\infty[\,,
\\ \nonumber
&&\ \ 
\sup_{x\in X}\sup_{r\in ]0,+\infty[}
\left\vert 
\int_{Y\setminus B(x,r)}K(x,y)\,d\nu(y)\right\vert<+\infty
 \biggl.\biggr\}  
\end{eqnarray*}
and
\begin{eqnarray*}
\lefteqn{
\|K\|_{{\mathcal{K}}_{ s_1, s_2, s_3  }^\sharp(X\times Y)}
\equiv
\|K\|_{{\mathcal{K}}_{ s_1, s_2, s_3  }(X\times Y)}
}
\\ \nonumber
&&+
\sup_{x\in X}\sup_{r\in ]0,+\infty[}
\left\vert
\int_{Y\setminus B(x,r)}K(x,y)\,d\nu(y)
\right\vert\quad\forall 
 K\in {\mathcal{K}}_{ s_1, s_2, s_3  }^\sharp(X\times Y)\,.
 \end{eqnarray*}
\end{defn}
Clearly,  $({\mathcal{K}}^\sharp_{ s_{1},s_{2},s_{3}   }(X\times Y),\|\cdot\|_{  {\mathcal{K}}^\sharp_{s_{1},s_{2},s_{3}   }(X\times Y)  })$ is a normed space and ${\mathcal{K}}^\sharp_{ s_{1},s_{2},s_{3}   }(X\times Y)$ is continuously embedded into ${\mathcal{K}}_{ s_{1},s_{2},s_{3}   }(X\times Y)$.  Then one can establish the H\"{o}lder continuity of $Q[Z,g,1]$ in the `singular' case $s_1=\upsilon_Y$ by means of the following statement, that extends some work of Gatto \cite[Proof of Thm.~3]{Gat09} (cf.~\cite[Prop.~6.3 (ii)]{La22a}).
\begin{prop}\label{prop:qzhq}
Let 
$
(M,d)$ be a metric space. Let $X$, $Y\subseteq M$. Let 
\[
\upsilon_Y\in ]0,+\infty[\,, \ \beta\in]0,1]\,,  \ 
 s_2\in  [\beta, +\infty[\,, \ s_3\in]0,1]\,.
\]
 Let $\nu$ be as in (\ref{eq:nu}), $\nu(Y)<+\infty$. Then the following statements hold.
\begin{enumerate}
\item[(b)] If $s_2-\beta>\upsilon_Y$, $s_2< \upsilon_Y+\beta+s_3$ and $Y$ is  upper $\upsilon_Y$-Ahlfors regular with respect to $X$, then the bilinear map
from 
\[
{\mathcal{K}}^\sharp_{ \upsilon_Y,s_{2},s_{3}   }(X\times Y)\times C^{0,\beta} (X\cup Y)
\quad\text{to}\quad
C_b^{0,\min\{
 \beta,
\upsilon_Y+s_3+\beta-s_2
\}}(X)\,,
\]
 which takes $(Z,g)$ to $Q[Z,g,1]$ is continuous. 
\item[(bb)] If $s_2-\beta=\upsilon_Y$ and $Y$ is strongly upper $\upsilon_Y$-Ahlfors regular with respect to $X$, then the bilinear map
from 
\[
{\mathcal{K}}^\sharp_{\upsilon_Y,s_{2},s_{3}   }(X\times Y)\times C^{0,\beta}(X\cup Y)
\quad\text{to}\quad C_b^{0,\max\{ r^\beta
,\omega_{s_3}(r)
\}}(X)\,,
\]
 which takes $(Z,g)$ to $Q[Z,g,1]$ is continuous. 
 
  \item[(bbb)] If $s_2-\beta<\upsilon_Y$ and $Y$ is  upper $\upsilon_Y$-Ahlfors regular with respect to $X$, then the bilinear map
from 
\[
{\mathcal{K}}^\sharp_{\upsilon_Y,s_{2},s_{3}   }(X\times Y)\times C^{0,\beta}(X\cup Y)
\quad\text{to}\quad C_b^{0,\min\{  \beta
, s_3 
\}}(X)\,,
\]
 which takes $(Z,g)$ to $Q[Z,g,1]$ is continuous. 
 \end{enumerate}
\end{prop}
Here $C^{0,\beta} (X\cup Y)$ denotes the semi-normed space of $\beta$-H\"{o}lder continuous functions on $X\cup Y$, $C^{0,\gamma}_b(X)$ denotes the normed space of bounded $\gamma$-H\"{o}lder continuous functions on $X$ for all $\gamma\in]0,1]$ and $C_b^{0,\max\{ r^\beta
,\omega_{s_3}(r)
\}}(X)$ denotes the normed space of bounded generalized H\"{o}lder continuous functions on $X$ with modulus of continuity $\max\{ r^\beta
,\omega_{s_3}(r)
\}$ (see notation around (\ref{om})--(\ref{omth}) below). Now the assumptions of Proposition \ref{prop:qzhq} require that the kernel $Z$ belongs to ${\mathcal{K}}^\sharp_{\upsilon_Y,s_{2},s_{3}   }(X\times Y)$ and the membership in ${\mathcal{K}}^\sharp_{\upsilon_Y,s_{2},s_{3}   }(X\times Y)$ requires an estimation of the maximal function of the kernel $Z$, \textit{i.e.}, to show that
\begin{equation}\label{eq:maxfubo}
\sup_{x\in X}\sup_{r\in ]0,+\infty[}
\left\vert
\int_{Y\setminus B(x,r)}Z(x,y)\,d\nu(y)
\right\vert <+\infty
\end{equation}
and here we wonder whether such assumption is actually necessary. Indeed, under the only assumption that the kernel $Z$ belongs to ${\mathcal{K}}_{\upsilon_Y,s_{2},s_{3}   }(X\times Y)$ and that $g$ is $\beta$-H\"{o}lder continuous, the integral that defines $Q[Z,g,1]$ is weakly singular.

Thus the main purpose of the present paper  is to establish what conditions on the maximal function of 
$Z$  as in (\ref{eq:maxfubo}) are actually necessary for the continuity of $Q[Z,\cdot,1]$ when $Z\in{\mathcal{K}}_{\upsilon_Y,s_{2},s_{3}   }(X\times Y)$. Namely, we prove the validity of 
  (\ref{eq:maxfuboa})--(\ref{eq:maxfuboc}), that are conditions on the growth of the left hand side of  (\ref{eq:maxfubo}) if $r$ is small (see Proposition \ref{prop:mfbd} (i)--(iii)). In order to prove such necessary conditions, we formulate an extra assumption on the set $X$, but no extra assumptions on the measure $\nu$ (see condition (\ref{prop:mfbd0})). 

We also note that one can formulate  some extra restrictions on the parameters in order that conditions (\ref{eq:maxfuboa})--(\ref{eq:maxfuboc}) imply the validity of  condition (\ref{eq:maxfubo}) and that the target space of $Q[Z,\cdot,1]$ equals $C_b^{0,\beta}(X)$.

To do so, in  case of statement (b),  one can   further require that  $s_2\leq \upsilon_Y+s_3$, an assumption   that is compatible with the condition $s_2=\upsilon_Y+s_3$ of standard kernels. 

Instead, in case of statement (bb), one can further require that $s_3>\beta$ and   in case of statement (bbb), one can further require that $s_3\geq\beta$, but such requirements are not compatible with the 
condition $s_2=\upsilon_Y+s_3$ of standard kernels.  

In all cases (b)---(bbb), we are saying   that  when   $Q[Z,\cdot,1]$ is bounded with no decrease of regularity, then the necessary condition (\ref{eq:maxfubo})  holds true. Instead when   $Q[Z,\cdot,1]$ is bounded with some decrease of regularity, then condition (\ref{eq:maxfubo}) is replaced by a growth condition of the maximal function in the variable $r$.

  In section \ref{sec:acodima}, we present an application of  Proposition \ref{prop:mfbd} (i) in  case $Y$ is a compact differentiable manifold.

\section{Analysis of the operator $Q$}
We first introduce  the following two technical lemmas, that generalize to case $X\neq Y$ the corresponding statements of Gatto \cite[p.~104]{Gat09}.
 For a proof, we refer to \cite[Lem.~3.2, 3.4, 3.6]{La22a}. \begin{lem}\label{lem:cominur}Let 
$
(M,d)$ be a metric space.  Let $X$, $Y\subseteq M$. 
Let $\upsilon_Y\in]0,+\infty[$. Let $\nu$ be as in (\ref{eq:nu}).    Let  $Y$ be  upper $\upsilon_Y$-Ahlfors regular with respect to $X$. Then the following statements hold.
\begin{enumerate}
\item[(i)] $\nu(\{x\})=0$ for all $x\in X\cap Y$.
\item[(ii)] Let $\nu(Y)<+\infty$. If $s\in ]0,\upsilon_Y[$, then
\begin{equation}\label{lem:cominur1}
c'_{s,X,Y}\equiv\sup_{x\in X}\int_Y\frac{d\nu(y)}{d(x,y)^s}
\leq \nu(Y)a^{-s}+c_{X,Y,\upsilon_Y}\frac{\upsilon_Y}{\upsilon_Y-s}a^{\upsilon_Y-s}
\end{equation}
for all $a\in]0,r_{X,Y,\upsilon_Y}[$.  If $s=0$, then 
\[
c'_{0,X,Y}\equiv\sup_{x\in X}\int_Y\frac{d\nu(y)}{d(x,y)^0}
= \nu(Y) \,.
\]	 
\item[(iii)]  Let  $\nu(Y)<+\infty$ whenever $r_{X,Y,\upsilon_Y}<+\infty$. If  $s\in]-\infty,\upsilon_Y[$, then
\[
c''_{s,X,Y}\equiv\sup_{(x,t)\in X\times]0,+\infty[}
t^{s-\upsilon_Y}\int_{B(x,t)\cap Y}\frac{d\nu(y)}{d(x,y)^s}<+\infty
\,.
\]
\end{enumerate}
\end{lem}
Then we have the following Lemma. 
\begin{lem}\label{lem:cominu} Let 
$
(M,d)$ be a metric space. Let $X$, $Y\subseteq M$. 
Let $\upsilon_Y\in]0,+\infty[$. Let $\nu$ be as in (\ref{eq:nu}), $\nu(Y)<+\infty$.  Then the following statements hold.
\begin{enumerate}
\item[(i)] Let  $Y$ be   upper $\upsilon_Y$-Ahlfors regular with respect to $X$.  If $s\in ]\upsilon_Y,+\infty[$, then 
\[
c'''_{s,X,Y}\equiv\sup_{(x,t)\in X\times]0,+\infty[}
t^{s-\upsilon_Y}\int_{Y\setminus B(x,t) }\frac{d\nu(y)}{d(x,y)^s}<+\infty\,.
\]
\item[(ii)] Let  $Y$ be strongly upper $\upsilon_Y$-Ahlfors regular with respect to $X$.   Then
\[
c^{iv}_{X,Y}\equiv\sup_{(x,t)\in X\times]0,1/e[}
\vert\log t\vert^{-1}\int_{Y\setminus B(x,t) }\frac{d\nu(y)}{d(x,y)^{\upsilon_Y}}<+\infty\,.
\]
\end{enumerate}
\end{lem}
Since we plan to analyze the operator $Q$ of (\ref{prop:qz6}) in H\"{o}lder spaces, we now introduce some notation. Let $X$ be a set. Then we set
\[
B(X)\equiv\left\{
f\in {\mathbb{C}}^X:\,f\ \text{is\ bounded}
\right\}
\,,\quad
\|f\|_{B(X)}\equiv\sup_X\vert f\vert \qquad\forall f\in B(X)\,.
\]
If $(M,d)$ is a metric space, then $C^0(M)$ denotes the set of continuous functions from $M$ to ${\mathbb{C}}$ and we introduce the subspace
$
C^0_b(M)\equiv C^0(M)\cap B(M)
$
of $B(M)$.  Let $\omega$ be a function from $[0,+\infty[$ to itself such that
\begin{eqnarray}
\nonumber
&&\qquad\qquad\omega(0)=0,\qquad \omega(r)>0\qquad\forall r\in]0,+\infty[\,,
\\
\label{om}
&&\qquad\qquad\omega\ {\text{is\   increasing,}}\ \lim_{r\to 0^{+}}\omega(r)=0\,,
\\
\nonumber
&&\qquad\qquad{\text{and}}\ \sup_{(a,t)\in[1,+\infty[\times]0,+\infty[}
\frac{\omega(at)}{a\omega(t)}<+\infty\,.
\end{eqnarray}
If $f$ is a function from a subset ${\mathbb{D}}$ of a metric space $(M,d)$   to ${\mathbb{C}}$,  then we denote by   $\vert f:{\mathbb{D}}\vert_{\omega (\cdot)}$  the $\omega(\cdot)$-H\"older constant  of $f$, which is delivered by the formula   
\[
\vert f:{\mathbb{D}}\vert_{\omega (\cdot)
}
\equiv
\sup\left\{
\frac{\vert f( x )-f( y)\vert }{\omega(d( x, y))
}: x, y\in {\mathbb{D}} ,  x\neq
 y\right\}\,.
\]        
If $\vert f:{\mathbb{D}}\vert_{\omega(\cdot)}<\infty$, we say that $f$ is $\omega(\cdot)$-H\"{o}lder continuous. Sometimes, we simply write $\vert f\vert_{\omega(\cdot)}$  
instead of $\vert f:{\mathbb{D}}\vert_{\omega(\cdot)}$. The
subset of $C^{0}({\mathbb{D}} ) $  whose
functions  are
$\omega(\cdot)$-H\"{o}lder continuous    is denoted  by  $C^{0,\omega(\cdot)} ({\mathbb{D}})$
and $\vert f:{\mathbb{D}}\vert_{\omega(\cdot)}$ is a semi-norm on $C^{0,\omega(\cdot)} ({\mathbb{D}})$.  
Then we consider the space  $C^{0,\omega(\cdot)}_{b}({\mathbb{D}} ) \equiv C^{0,\omega(\cdot)} ({\mathbb{D}} )\cap B({\mathbb{D}} ) $ with the norm \[
\|f\|_{ C^{0,\omega(\cdot)}_{b}({\mathbb{D}} ) }\equiv \sup_{x\in {\mathbb{D}} }\vert f(x)\vert
+\vert f\vert_{\omega(\cdot)}\qquad\forall f\in C^{0,\omega(\cdot)}_{b}({\mathbb{D}} )\,.
\] 
In the case in which $\omega(\cdot)$ is the function 
$r^{\alpha}$ for some fixed $\alpha\in]0,1]$, a so-called H\"{o}lder exponent, we simply write $\vert \cdot:{\mathbb{D}}\vert_{\alpha}$ instead of
$\vert \cdot:{\mathbb{D}}\vert_{r^{\alpha}}$, $C^{0,\alpha} ({\mathbb{D}})$ instead of $C^{0,r^{\alpha}} ({\mathbb{D}})$, $C^{0,\alpha}_{b}({\mathbb{D}})$ instead of $C^{0,r^{\alpha}}_{b} ({\mathbb{D}})$, and we say that $f$ is $\alpha$-H\"{o}lder continuous provided that 
$\vert f:{\mathbb{D}}\vert_{\alpha}<\infty$.
Next we introduce a function that we need for a generalized H\"{o}lder norm. For each $\theta\in]0,1]$, we define the function $\omega_{\theta}(\cdot)$ from $[0,+\infty[$ to itself by setting
\begin{equation}\label{omth}
\omega_{\theta}(r)\equiv
\left\{
\begin{array}{ll}
0 &r=0\,,
\\
r^{\theta}\vert \ln r \vert  &r\in]0,r_{\theta}]\,,
\\
r_{\theta}^{\theta}\vert \ln r_{\theta} \vert  & r\in ]r_{\theta},+\infty[\,,
\end{array}
\right.
\end{equation}
where
$
r_{\theta}\equiv e^{-1/\theta}
$ for all $\theta\in ]0,1]$. Obviously, $\omega_{\theta}(\cdot) $ is concave and satisfies   condition (\ref{om}).
 We also note that if ${\mathbb{D}}\subseteq M$, then the continuous embeddings
\[
C^{0, \theta }_b({\mathbb{D}})\subseteq 
C^{0,\omega_\theta(\cdot)}_b({\mathbb{D}})\subseteq 
C^{0,\theta'}_b({\mathbb{D}})
\]
hold  for all $\theta'\in ]0,\theta[$. We now turn to analyze the integral operator $Q$ of (\ref{prop:qz6}). As we have already said, the sufficient conditions of Proposition \ref{prop:qzhq}
   require that the kernel $Z$ belongs to ${\mathcal{K}}^\sharp_{\upsilon_Y,s_{2},s_{3}   }(X\times Y)$ and the membership in ${\mathcal{K}}^\sharp_{\upsilon_Y,s_{2},s_{3}   }(X\times Y)$ requires an estimation of the maximal function of the kernel $Z$, \textit{i.e.}, to show the validity of condition (\ref{eq:maxfubo}) and we now prove the following `converse' statement that shows that under certain assumptions condition (\ref{eq:maxfubo}) is actually necessary.
\begin{prop}\label{prop:mfbd}
Let 
$
(M,d)$ be a metric space.  Let $X$, $Y\subseteq M$. Assume that 
\begin{eqnarray}\nonumber
&&\text{there\ exists}\  a\in]0,+\infty[\ \text{such\ that\  for\ each}\  \rho\in]0,a[\ \text{and}\  x'\in X,
 \\ \label{prop:mfbd0}
&&\text{there\ exists}\  x''\in X\ \text{such\ that}\  d(x',x'')=\rho\,.
\end{eqnarray}
 Let 
\[
\upsilon_Y\in ]0,+\infty[\,, \ \beta\in]0,1]\,,  \ 
 s_2\in  [\beta, +\infty[\,, \ s_3\in]0,1]\,.
\]
 Let $\nu$ be as in (\ref{eq:nu}), $\nu(Y)<+\infty$. Let $Z\in {\mathcal{K}}_{\upsilon_Y,s_{2},s_{3}   }(X\times Y)$. Then
\begin{equation}\label{prop:mfbd1}
Z(x,\cdot)\ \text{is}\ \nu-\text{integrable\ in}\ Y\setminus B(x,r)
\ \text{for\ all}\ (x,r)\in X\times]0,+\infty[ 
\end{equation}
and the following statements hold.
\begin{enumerate}
\item[(i)] Let  $s_2-\beta>\upsilon_Y$, $s_2< \upsilon_Y+\beta+s_3$. Let $Y$ be  upper $\upsilon_Y$-Ahlfors regular with respect to $X$. 
If the linear map
from $C^{0,\beta} (X\cup Y)$ to $C_b^{0,\min\{
 \beta,
\upsilon_Y+s_3+\beta-s_2
\}}(X)$ that takes $g$ to  $Q[Z,g,1]$ is continuous, then 
\begin{equation}\label{eq:maxfuboa}
\sup_{x\in X}\sup_{r\in ]0,e^{-1/s_3}[}
\left\vert
\int_{Y\setminus B(x,r)}Z(x,y)\,d\nu(y)
\right\vert  \frac{r^\beta}{\max\{r^\beta,r^{\upsilon_Y+s_3+\beta-s_2}
\}}
<+\infty\,.
\end{equation}
In particular, if $s_2\leq \upsilon_Y+s_3$, then $\min\{
 \beta,
\upsilon_Y+s_3+\beta-s_2
\}=\beta$ and inequality (\ref{eq:maxfubo}) holds true.
\item[(ii)] Let $s_2-\beta=\upsilon_Y$. Let $Y$ be strongly upper $\upsilon_Y$-Ahlfors regular with respect to $X$. If the linear map
from $C^{0,\beta} (X\cup Y)$ to $C_b^{0,\max\{ r^\beta
,\omega_{s_3}(r)
\}}(X)$ that takes $g$ to  $Q[Z,g,1]$ is continuous, then 
\begin{equation}\label{eq:maxfubob}
\sup_{x\in X}\sup_{r\in ]0,e^{-1/s_3}[}
\left\vert
\int_{Y\setminus B(x,r)}Z(x,y)\,d\nu(y)
\right\vert  \frac{r^\beta}{\max\{ r^\beta
,\omega_{s_3}(r)
\}}
<+\infty\,.
\end{equation} 
In particular, if  $s_3>\beta$, then $C_b^{0,\max\{ r^\beta
,\omega_{s_3}(r)
\}}(X)=C_b^{0,\beta} (X)$ and  inequality (\ref{eq:maxfubo}) holds true.
\item[(iii)] Let $s_2-\beta<\upsilon_Y$. Let $Y$ be  upper $\upsilon_Y$-Ahlfors regular with respect to $X$.  If the linear map
from $C^{0,\beta} (X\cup Y)$ to $C_b^{0,\min\{  \beta
, s_3 
\}}(X)$  that takes $g$ to  $Q[Z,g,1]$ is continuous, then 
\begin{equation}\label{eq:maxfuboc}
\sup_{x\in X}\sup_{r\in ]0,e^{-1/s_3}[}
\left\vert
\int_{Y\setminus B(x,r)}Z(x,y)\,d\nu(y)
\right\vert  \frac{r^\beta}{\max\{ r^\beta
,r^{s_3}
\}}
<+\infty\,.
\end{equation} 
In particular, if  $s_3\geq \beta$, then 
$C_b^{0,\min\{  \beta
, s_3 
\}}(X)=C^{0,\beta}_b (X)$ and 
inequality (\ref{eq:maxfubo}) holds true.
\end{enumerate}
\end{prop}
{\bf Proof.}  Since $Z\in {\mathcal{K}}_{\upsilon_Y,s_{2},s_{3}   }(X\times Y)$, we have
\begin{equation}\label{prop:mfbd2}
\vert Z(x,y)\vert \leq\frac{\|Z\|_{{\mathcal{K}}_{\upsilon_Y,s_{2},s_{3}   }(X\times Y)}}{d(x,y)^{\upsilon_Y}}
\leq
\frac{\|Z\|_{{\mathcal{K}}_{\upsilon_Y,s_{2},s_{3}   }(X\times Y)}}{r^{\upsilon_Y}}\qquad\forall y\in Y\setminus B(x,r)
\end{equation}
for all $(x,r)\in X\times]0,+\infty[$. On the other hand,   $Z$ is continuous and $\nu(Y)$ is finite and thus condition (\ref{prop:mfbd1}) is satisfied.
Inequality (\ref{prop:mfbd2}) also shows that it suffices to show inequalities
(\ref{eq:maxfubo}), (\ref{eq:maxfuboa})--(\ref{eq:maxfuboc}),  for $r\in]0, \min\{a,e^{-1/s_3}\}/2[$ ($\subseteq]0,\min\{a,e^{-1}\}/2[$). We find convenient to set
\[
\phi(r)\equiv\left\{
\begin{array}{ll}
 \max\{
 r^\beta,
r^{\upsilon_Y+s_3+\beta-s_2}
\} &\text{if}\  s_2-\beta>\upsilon_Y,  s_2< \upsilon_Y+\beta+s_3\,,
\\
\max\{ r^\beta
,\omega_{s_3}(r)
\}
&\text{if}\  s_2-\beta=\upsilon_Y\,,
\\
\max\{ r^\beta
,r^{s_3}
\}&\text{if}\  s_2-\beta<\upsilon_Y
\end{array}
\right.
\]
for all $r\in[0,+\infty[$. The idea consists in proving an upper estimate for
\[
\left\vert
\int_{Y\setminus B(x',2d(x',x''))}Z(x',y)\,d\nu(y)
\right\vert  \vert g(x')-g(x'')\vert   
\]
in terms of $\phi (d(x',x''))$, of the norm of the kernel $Z$ and of the seminorm of $g$   
when $g\in C^{0,\beta} (X\cup Y)$, $x'$, $x''\in X$, $d(x',x'')\leq \min\{a,e^{-1/s_3}\}/2$ and then in making a suitable choice of $g$ that enables to establish an upper bound for $\left\vert
\int_{Y\setminus B(x',2d(x',x''))}Z(x',y)\,d\nu(y)
\right\vert $ that is independent of $x'$ and $x''$ and then to choose $x''$ so that $2d(x',x'')=r$ and to obtain the upper bound we need for $\left\vert
\int_{Y\setminus B(x',r)}Z(x',y)\,d\nu(y)
\right\vert $. To do so, we first look at the following elementary equality
\begin{eqnarray*}
\lefteqn{
Q[Z,g,1](x')-Q[Z,g,1](x'')
}
\\ \nonumber
&&\qquad
=\int_YZ(x',y)(g(x')-g(y))\,d\nu(y)
-
\int_YZ(x'',y)(g(x'')-g(y))\,d\nu(y)
\\ \nonumber
&&\qquad
=\int_{Y\cap B(x',2d(x',x''))}
 Z(x',y)(g(x')-g(y))\,d\nu(y)
 \\ \nonumber
&&\qquad\quad
-\int_{Y\cap B(x',2d(x',x''))}Z(x'',y)(g(x'')-g(y))\,d\nu(y)
\\ \nonumber
&&\qquad\quad
+\int_{Y\setminus B(x',2d(x',x''))}Z(x',y)
[(g(x')-g(y)) -(g(x'')-g(y)) ]\,d\nu(y)
\\ \nonumber
&&\qquad\quad
+\int_{Y\setminus B(x',2d(x',x''))}[Z(x',y)
-Z(x'',y)](g(x'')-g(y))\,d\nu(y) 
\end{eqnarray*}
for all $x',x''\in X$, $d(x',x'')\leq \min\{a,e^{-1/s_3}\}/2$ and $g\in C^{0,\beta} (X\cup Y)$. Next we note that under each of the assumptions of statements (i)--(iii), the linear operator
$Q[Z,\cdot,1]$ is continuous from $C^{0,\beta} (X\cup Y)$ to 
$C^{0,\phi(\cdot)}_b (X)$. Then we find convenient to introduce the symbol 
${\mathcal{L}}(C^{0,\beta} (X\cup Y),C^{0,\phi(\cdot)}_b (X))$ for the normed space of linear and continuous operators from $C^{0,\beta} (X\cup Y)$ to $C^{0,\phi(\cdot)}_b (X)$.  Then Lemma \ref{lem:cominur}   and the inclusion
$
B(x',2d(x',x''))\subseteq B(x'',3d(x',x''))
$
imply that 
\begin{eqnarray}\label{prop:mfbd2a}
\lefteqn{
\left\vert\int_{Y\setminus B(x',2d(x',x''))}Z(x',y)
[g(x'')- g(x')  ]\,d\nu(y)\right\vert 
}
\\ \nonumber
&&
\leq \|Q[Z,\cdot,1]\|_{{\mathcal{L}}(C^{0,\beta} (X\cup Y),C^{0,\phi(\cdot)}_b (X))}
\vert g:\,X\cup Y\vert_\beta \phi(d(x',x''))
\\ \nonumber
&&\quad
+\int_{Y\cap B(x',2d(x',x''))}
 \vert Z(x',y)\vert \,\vert g(x')-g(y)\vert \,d\nu(y)
 \\ \nonumber
&&\quad
+\int_{Y\cap B(x',2d(x',x''))}\vert Z(x'',y)\vert \,\vert g(x'')-g(y)\vert \,d\nu(y)
\\ \nonumber
&&\quad
 +\|Z\|_{{\mathcal{K}}_{\upsilon_Y,s_2,s_3}(X\times Y) }
\vert g:X\cup Y\vert_\beta 
\\ \nonumber
&&\qquad\quad
\int_{Y\setminus B(x',2d(x',x''))}
\frac{d(x',x'')^{s_3}}{d(x',y)^{s_2}}d(y,x'')^\beta\,d\nu(y)
\\ \nonumber
&&
\leq \|Q[Z,\cdot,1]\|_{{\mathcal{L}}(C^{0,\beta} (X\cup Y),C^{0,\phi(\cdot)}_b (X))}
\vert g:\,X\cup Y\vert_\beta \phi(d(x',x''))
\\ \nonumber
&&\quad
+\|Z\|_{{\mathcal{K}}_{\upsilon_Y,s_2,s_3}(X\times Y) }
\vert g:X\cup Y\vert_\beta 
\biggl\{\biggr.
\int_{ Y\cap B(x',2d(x',x'')) }
\frac{d\nu(y)}{d(x',y)^{\upsilon_Y-\beta}}
\\ \nonumber
&&\quad
+
\int_{ Y\cap B(x'',3d(x',x'')) }
\frac{d\nu(y)}{d(x'',y)^{\upsilon_Y-\beta}}
\biggl.\biggr\}
\\ \nonumber
&&\quad
 +\|Z\|_{{\mathcal{K}}_{\upsilon_Y,s_2,s_3}(X\times Y) }
\vert g:X\cup Y\vert_\beta 
\\ \nonumber
&&\qquad\quad
\times\int_{Y\setminus B(x',2d(x',x''))}
\frac{d(x',x'')^{s_3}}{d(x',y)^{s_2}}d(y,x'')^\beta\,d\nu(y)
\\ \nonumber
&&
\leq \|Q[Z,\cdot,1]\|_{{\mathcal{L}}(C^{0,\beta} (X\cup Y),C^{0,\phi(\cdot)}_b (X))}
\vert g:\,X\cup Y\vert_\beta \phi(d(x',x''))
\\ \nonumber
&&\quad
+\|Z\|_{{\mathcal{K}}_{\upsilon_Y,s_2,s_3}(X\times Y) }
\vert g:X\cup Y\vert_\beta c_{\upsilon_Y-\beta,X,Y}''
\\ \nonumber
&&\qquad\qquad
\times
\biggl\{\biggr.
 (2d(x',x''))^{\upsilon_Y-(\upsilon_Y-\beta) }
+
  (3d(x',x''))^{\upsilon_Y-(\upsilon_Y-\beta) } 
\biggl.\biggr\}
\\ \nonumber
&&\quad
 +\|Z\|_{{\mathcal{K}}_{\upsilon_Y,s_2,s_3}(X\times Y) }
\vert g:X\cup Y\vert_\beta 
\\ \nonumber
&&\qquad\quad
\times \int_{Y\setminus B(x',2d(x',x''))}
\frac{d(x',x'')^{s_3}}{d(x',y)^{s_2}}d(y,x'')^\beta\,d\nu(y)
\end{eqnarray}
for all $x',x''\in X$, $d(x',x'')\leq \min\{a,e^{-1/s_3}\}/2$ and $g\in C^{0,\beta} (X\cup Y)$. Next we note that the triangular inequality 
implies that
\[
d(y,x'')\leq d(y,x')+d(x',x'')
\leq 2d(x',y)\qquad \forall y\in Y\setminus B(x',2d(x',x'')) 
\]
and that accordingly
\begin{eqnarray}\label{prop:qzhq2}
\lefteqn{
\int_{Y\setminus B(x',2d(x',x'')}
\frac{d(x',x'')^{s_3}}{d(x',y)^{s_2}}d(y,x'')^\beta\,d\nu(y)
}
\\ \nonumber
&&\qquad
\leq 2^\beta
\int_{Y\setminus B(x',2d(x',x'')}
\frac{d\nu(y)}{d(x',y)^{s_2-\beta} }d(x',x'')^{s_3} 
\end{eqnarray}
for all $x',x''\in X$, $d(x',x'')\leq \min\{a,e^{-1/s_3}\}/2$. We now distinguish three cases. If $s_2-\beta>\upsilon_Y$ as in statement (i), then Lemma \ref{lem:cominu} (i)
   implies that
\begin{eqnarray}\label{prop:qzhq3}
\lefteqn{
\int_{Y\setminus B(x',2d(x',x''))}
\frac{d\nu(y)}{ d(x',y)^{s_2-\beta} }d(x',x'')^{s_3}
}
\\ \nonumber
&&\qquad
\leq c'''_{s_2-\beta,X,Y}(2d(x',x''))^{\upsilon_Y-(s_2-\beta)}d(x',x'')^{s_3}
\\ \nonumber
&&\qquad
\leq c'''_{
s_2-\beta,X,Y}d(x',x'')^{\upsilon_Y+s_3+\beta-s_2 } 
\end{eqnarray}
for all $x',x''\in X$, $d(x',x'')\leq \min\{a,e^{-1/s_3}\}/2$. If $s_2-\beta=\upsilon_Y$  as in statement (ii), then Lemma \ref{lem:cominu} (ii)   implies that
\begin{eqnarray}\label{prop:qzhq4}
\lefteqn{
\int_{Y\setminus B(x',2d(x',x''))}
\frac{d\nu(y)}{ d(x',y)^{s_2-\beta} }d(x',x'')^{s_3}
}
\\ \nonumber
&&\qquad
\leq c^{iv}_{X,Y}\vert \log (2d(x',x''))\vert \,d(x',x'')^{s_3}
\\ \nonumber
&&\qquad
\leq c^{iv}_{X,Y}\vert \log d(x',x'') \vert \,d(x',x'')^{s_3}
\left(
1+\frac{\log 2}{\vert \log d(x',x'') \vert }
\right)
\\ \nonumber
&&\qquad
\leq
2c^{iv}_{X,Y}\vert \log d(x',x'') \vert \,d(x',x'')^{s_3} 
\end{eqnarray}
for all $x',x''\in X$, $d(x',x'')\leq \min\{a,e^{-1/s_3}\}/2$.
If $s_2-\beta<\upsilon_Y$  as in statement (iii), then Lemma \ref{lem:cominur} (i)   implies that
\begin{equation}\label{prop:qzhq4a}
\int_{Y\setminus B(x',2d(x',x''))}
\frac{d\nu(y)}{ d(x',y)^{s_2-\beta} }d(x',x'')^{s_3}
\leq
c'_{s_2-\beta,X,Y}d(x',x'')^{s_3} 
\end{equation}
for all $x',x''\in X$, $d(x',x'')\leq \min\{a,e^{-1/s_3}\}/2$.

Under the assumptions of statement (i), we have 
$
s_2-\beta>\upsilon_Y 
$
 and inequalities (\ref{prop:mfbd2a}), (\ref{prop:qzhq2}), (\ref{prop:qzhq3}) imply that there exist $c_{(i)}\in]0,+\infty[$ such that
\begin{eqnarray*}
\lefteqn{
\left\vert\int_{Y\setminus B(x',2d(x',x''))}Z(x',y)
[g(x'')- g(x')  ]\,d\nu(y)\right\vert 
}
\\ \nonumber
&&\qquad
\leq c_{(i)} 
\vert g:X\cup Y\vert_\beta
\\ \nonumber
&&\qquad\quad
\times
\max\left\{\phi(d(x',x'')),
d(x',x'')^\beta,d(x',x'')^{\upsilon_Y+s_3+\beta-s_2}
\right\}
\\ \nonumber
&&\qquad
=c_{(i)} 
\vert g:X\cup Y\vert_\beta \phi(d(x',x'')),
\end{eqnarray*}
for all $x',x''\in X$, $d(x',x'')\leq \min\{a,e^{-1/s_3}\}/2$ and $g\in C^{0,\beta} (X\cup Y)$. 

Under the assumptions of statement (ii), we have $s_2-\beta=\upsilon_Y$   and inequalities (\ref{prop:mfbd2a}), (\ref{prop:qzhq2}), (\ref{prop:qzhq4}) imply that there exist $c_{(ii)}\in]0,+\infty[$ such that
\begin{eqnarray*}
\lefteqn{
\left\vert\int_{Y\setminus B(x',2d(x',x''))}Z(x',y)
[g(x'')- g(x')  ]\,d\nu(y)\right\vert 
}
\\ \nonumber
&&\qquad
\leq c_{(ii)} 
\vert g:X\cup Y\vert_\beta
\\ \nonumber
&&\qquad\quad
\times
\max\left\{
\phi(d(x',x'')),d(x',x'')^\beta,\vert \log d(x',x'') \vert \,d(x',x'')^{s_3}
\right\}
\\ \nonumber
&&\qquad
=c_{(ii)} 
\vert g:X\cup Y\vert_\beta \phi(d(x',x''))
\end{eqnarray*}
for all $x',x''\in X$, $d(x',x'')\leq \min\{a,e^{-1/s_3}\}/2$ and $g\in C^{0,\beta} (X\cup Y)$. 

Under the assumptions of statement (iii), we have $s_2-\beta<\upsilon_Y$ and inequalities (\ref{prop:mfbd2a}), (\ref{prop:qzhq2}), (\ref{prop:qzhq4a}) imply that there exist $c_{(iii)}\in]0,+\infty[$ such that
\begin{eqnarray*}
\lefteqn{
\left\vert\int_{Y\setminus B(x',2d(x',x''))}Z(x',y)
[g(x'')- g(x')  ]\,d\nu(y)\right\vert 
}
\\ \nonumber
&&\qquad
\leq c_{(iii)} 
\vert g:X\cup Y\vert_\beta
\\ \nonumber
&&\qquad\quad
\times
\max\left\{\phi(d(x',x'')), 
d(x',x'')^\beta, d(x',x'')^{s_3}
\right\}
\\ \nonumber
&&\qquad
=c_{(iii)} 
\vert g:X\cup Y\vert_\beta \phi(d(x',x''))
\end{eqnarray*}
for all $x',x''\in X$, $d(x',x'')\leq \min\{a,e^{-1/s_3}\}/2$ and $g\in C^{0,\beta} (X\cup Y)$. 

Hence, in all cases (i)--(iii), there exists $c\in]0,+\infty[$ such that
\begin{eqnarray}\label{prop:mfbd5}
\lefteqn{
\left\vert\int_{Y\setminus B(x',2d(x',x''))}Z(x',y)
 \,d\nu(y)\right\vert \,\vert g(x'')- g(x')\vert 
}
\\ \nonumber
&&\qquad
\leq c  
\vert g:X\cup Y\vert_\beta \phi(d(x',x''))
\end{eqnarray}
for all $x',x''\in X$, $d(x',x'')\leq \min\{a,e^{-1/s_3}\}/2$ and $g\in C^{0,\beta} (X\cup Y)$. Next we make a choice of $g$ that depends upon $x'$. Namely, we set
\[
g_{x'}(y)\equiv d(x',y)^\beta\qquad\forall y\in X\cup Y 
\]
for all $x'\in X$. By the elementary inequality
\[
\vert g_{x'}(u_1)-g_{x'}(u_2)\vert =\vert d(x',u_1)^\beta-d(x',u_2)^\beta\vert \leq 
d(u_1,u_2)^\beta\qquad\forall u_1, u_2\in X\cup Y\,,
\]
that follows by the Yensen inequality, we have
\[
\vert g_{x'}:X\cup Y\vert_\beta\leq 1\qquad\forall x'\in X\,.
\]
(cf., \textit{e.g.}, Folland \cite[Prop. 6.11]{Fo99}). Then $g_{x'}\in C^{0,\beta} (X\cup Y)$ and inequality (\ref{prop:mfbd5}) implies that 
\begin{eqnarray}\label{prop:mfbd6}
\lefteqn{
\left\vert\int_{Y\setminus B(x',2d(x',x''))}Z(x',y)
 \,d\nu(y)\right\vert d(x',x'')^\beta
 }
\\ \nonumber
&&\qquad
 =
\left\vert\int_{Y\setminus B(x',2d(x',x''))}Z(x',y)
 \,d\nu(y)\right\vert \,\vert g_{x'}(x'')- g_{x'}(x')\vert 
\\ \nonumber
&&\qquad
\leq c  
\vert g_{x'}:X\cup Y\vert_\beta \phi(d(x',x''))
\leq c \phi(d(x',x''))
\end{eqnarray}
for all $x',x''\in X$, $d(x',x'')\leq \min\{a,e^{-1/s_3}\}/2$. Thus if $(x',r)$ belongs to $ X\times]0,\min\{a,e^{-1/s_3}\}/2[$ our assumption implies that there exists $x''\in X$ such that $d(x',x'')=r/2$ and thus inequality (\ref{prop:mfbd5}) implies that
\[
\left\vert\int_{Y\setminus B(x',r)}Z(x',y)
 \,d\nu(y)\right\vert  r^\beta 
  2^{-\beta}	 
 \leq c  \phi(r/2)\leq  c  \phi(r)
\]
and thus the proof is complete. The last part of statements (i)--(iii) is an immediate consequence of the corresponding first parts.  \hfill  $\Box$ 

\vspace{\baselineskip}

We note that the condition of existence of $a$ implies that each point of $X$ is an accumulation point for $X$.

\section{An application in  case $Y$ is a compact differentiable manifold}
\label{sec:acodima}
Since a compact manifold $Y$ of class $C^1$ that is imbedded in ${\mathbb{R}}^n$   is    $(n-1)$-upper Ahlfors regular, we now present an application of Proposition \ref{prop:mfbd} (i) in case $Y$ is a compact manifold  of class $C^1$ that is imbedded in ${\mathbb{R}}^n$. To do so, we need the following elementary lemma, that shows that $Y$ satisfies the technical condition (\ref{prop:mfbd0}).
\begin{lem}\label{lem:tecolo}
 Let $n\in {\mathbb{N}}\setminus\{0\}$. Let $Y$ be a compact manifold of class $C^0$ that is imbedded in ${\mathbb{R}}^n$ and of dimension $m$. Let $m\geq 1$. 
 Then $Y$ satisfies condition (\ref{prop:mfbd0}) with $X=Y$.
\end{lem}
{\bf Proof.}  Since $Y$ is a compact manifold of class $C^0$, $Y$ can be covered by a finite number of open connected domains of charts, each of which cannot be equal to $Y$.  Then by taking $a$ to be one half of a Lebesgue number for such a finite open cover, for each $x'\in Y$ and $\rho\in]0,a[$, 
the set $Y\cap\overline{{\mathbb{B}}_n(x',\rho)} $ is contained in at least one open connected domain of chart of the finite cover of $Y$, say $A$  (cf.~\textit{e.g.}, Dugundji~\cite[Theorem~4.5, Chap.~XI]{Du76}). Since $A$ 
is homeomorphic to a open subset of ${\mathbb{R}}^m$ that is not empty, $A$ cannot be compact. Since  $Y\cap\overline{{\mathbb{B}}_n(x',\rho)} $ is    compact, $Y\cap\overline{{\mathbb{B}}_n(x',\rho)} $ cannot be equal to $A$ and thus the  set $A\setminus (Y\cap\overline{{\mathbb{B}}_n(x',\rho)})$ cannot be empty. Since the  
$Y\cap {\mathbb{B}}_n(x',\rho)$ contains $x'$ it is not empty.
Since the  
$Y\cap {\mathbb{B}}_n(x',\rho)$ is open in $A$ and is not empty and $A$ is connected, we conclude that  $Y\cap\partial{\mathbb{B}}_n(x',\rho)$ cannot be empty and condition (\ref{prop:mfbd0}) with $X=Y$ holds true.\hfill  $\Box$ 

\vspace{\baselineskip}

Next we introduce the following. For the definition of tangential gradient $ {\mathrm{grad}}_{Y}$, we refer \textit{e.g.},   to Kirsch and Hettlich \cite[A.5]{KiHe15}, Chavel~\cite[Chap.~1]{Cha84}. 
 For a proof, we refer to \cite[Thm.~6.2]{La22b}
\begin{thm}\label{thm:tgkwmcb}
Let $n\in {\mathbb{N}}$,  $n\geq 2$. Let $Y$ be a compact manifold of class $C^1$ that is imbedded in ${\mathbb{R}}^n$. Let $s_1\in [0,(n-1)[$. Let $\beta\in ]0,1]$, $t_1\in ]0,(n-1)+\beta[$. Let the kernel $K\in {\mathcal{K}}_{s_1,s_1+1,1}(Y\times Y)$ satisfy the following assumptions
\begin{eqnarray*}
 &&K(\cdot,y)\in C^1(Y\setminus\{y\})  \qquad\forall y\in Y\,,
\qquad
 \int_YK(\cdot,y)\,d\sigma_y\in C^1(Y) \,,
 \\
 &&{\mathrm{grad}}_{Y,x}K(\cdot,\cdot)\in\left(
 {\mathcal{K}}_{t_1,Y\times Y} 
 \right)^n\,.
\end{eqnarray*}
Let $\mu\in C^{0,\beta}_b(Y)$. Then the function $\int_YK(\cdot,y)\mu(y)\,d\sigma_y$ is of class $C^1(Y) $, the function $ [{\mathrm{grad}}_{Y,x}  K(x,y)](\mu(y)-\mu(x))$ is  integrable in  $y\in Y$ for all $x\in Y$ and    formula 
\begin{eqnarray}\label{thm:tadelg1}
\lefteqn{
{\mathrm{grad}}_{Y,x} \int_YK(x,y)\mu(y)\,d\sigma_y
}
\\ \nonumber
&&\ \ 
=\int_Y[{\mathrm{grad}}_{Y,x}K(x,y)](\mu(y)-\mu(x))\,d\sigma_y+\mu(x){\mathrm{grad}}_{Y} \int_YK(x,y) \,d\sigma_y \,,
\end{eqnarray}
for all $x\in Y$,   for the tangential gradient of $\int_YK(\cdot,y)\mu(y)\,d\sigma_y$ holds true.
\end{thm}
By combining Proposition \ref{prop:qzhq} (b) and Theorem \ref{thm:tgkwmcb}, one can  prove the following  continuity Theorem \ref{thm:iokreg} for the integral operator with kernel $K$ and with values into a Schauder space on a compact manifold  $Y$ of class $C^1$ (cf.~\cite[Thm.~6.3 (ii) (b)]{La22b}). For the definition of the Schauder space
$C^{1,\beta}(Y)$ of functions $\mu$ of class $C^1$ on $Y$ such that the tangential gradient of $\mu$ is $\beta$-H\"{o}lder continuous or for an equivalent definition based on a finite family of parametrizations of $Y$, we refer for example to \cite[\S 2.20]{DaLaMu21}.
 \begin{thm}\label{thm:iokreg}
 Let $n\in {\mathbb{N}}$,  $n\geq 2$.  Let $Y$ be a compact manifold of class $C^1$ that is imbedded in ${\mathbb{R}}^n$.  Let $s_1\in [0,(n-1)[$. 
 Let $\beta\in ]0,1]$, $t_1\in[\beta,(n-1)+\beta[$, $t_2\in [ \beta,+\infty[$, $t_3\in]0,1]$. Let the kernel $K\in {\mathcal{K}}_{s_1,s_1+1,1}(Y\times Y)$ satisfy the following assumption 
 \[
 K(\cdot,y)\in C^1(Y\setminus\{y\})  \quad\forall y\in Y\,.
 \]
Let $t_1=(n-1)$, ${\mathrm{grad}}_{Y,x}K(\cdot,\cdot)\in\left(
 {\mathcal{K}}^\sharp_{t_1,t_2,t_3}(Y\times Y)
 \right)^n$,  $t_2-\beta>(n-1)$, $t_2<(n-1)+\beta+t_3$ and 
\[
 \int_YK(\cdot,y)\,d\nu(y)\in C^{1,\min\{  \beta, (n-1)+t_3+\beta-t_2\} }(Y)\,.
\]
Then the map from $C^{0,\beta}_b(Y)$ to $C^{1,\min\{  \beta, (n-1)+t_3+\beta-t_2\}}_b(Y)$ that takes $\mu$ to the function  $\int_YK(\cdot,y)\mu(y)\,d\sigma_y$ is linear and continuous. 
\end{thm}
By combining Proposition \ref{prop:mfbd} (i) and Theorem \ref{thm:tgkwmcb}, we can now prove the following  converse result of the continuity Theorem \ref{thm:iokreg} with the additional restriction that $t_2\leq(n-1)+t_3$, a case in which $C^{1,\min\{  \beta, (n-1)+t_3+\beta-t_2\} }(Y)=C^{1,\beta }(Y)$.
We also set
\[
 {\mathbb{B}}_{n}(x,\rho)\equiv \left\{
y\in {\mathbb{R}}^{n}:\,\vert x-y\vert <\rho
\right\}   
\]
For all $\rho>0$, $ x\in{\mathbb{R}}^{n}$.
\begin{thm}\label{thm:iokregn}
 Let $n\in {\mathbb{N}}$,  $n\geq 2$. Let $Y$ be a compact manifold of class $C^1$ that is imbedded in ${\mathbb{R}}^n$.  Let $s_1\in [0,(n-1)[$. 
 Let $\beta\in ]0,1]$, $t_1\in[\beta,(n-1)+\beta[$, $t_2\in [ \beta,+\infty[$, $t_3\in]0,1]$. Let the kernel $K\in {\mathcal{K}}_{s_1,s_1+1,1}(Y\times Y)$ satisfy the following assumption 
 \[
 K(\cdot,y)\in C^1(Y\setminus\{y\})  \quad\forall y\in Y\,.
 \]
Let $t_1=(n-1)$, ${\mathrm{grad}}_{Y,x}K(\cdot,\cdot)\in\left(
 {\mathcal{K}}_{t_1,t_2,t_3}(Y\times Y)
 \right)^n$,   $t_2-\beta>(n-1)$, $t_2\leq(n-1)+t_3$. 
If the map from $C^{0,\beta}(Y)$ to $C^{1,\beta}(Y)$ that takes $\mu$ to the function  $\int_YK(\cdot,y)\mu(y)\,d\sigma_y$ is linear and continuous, then inequality
\begin{equation}\label{thm:iokregn1}
\sup_{x\in Y}\sup_{r\in ]0,+\infty[}
\left\vert
\int_{Y\setminus {\mathbb{B}}_n(x,r)}{\mathrm{grad}}_{Y,x}K(x,y)\,d\sigma_y
\right\vert <+\infty
\end{equation}
holds true.
\end{thm}
{\bf Proof.} One can easily verify that $Y$ is upper $(n-1)$-Ahlfors regular with respect to $Y$. By Lemma \ref{lem:tecolo}, $Y$ satisfies condition (\ref{prop:mfbd0}). Since $1\in C^{0,\beta}(Y)$, our assumption implies that 
\[
 \int_YK(\cdot,y)\,d\sigma_y\in C^{1,\beta}(Y)\,.
\]
By Theorem \ref{thm:tgkwmcb}, the function $\int_YK(\cdot,y)\mu(y)\,d\sigma_y$ is of class $C^1(Y)$ for all $\mu \in C^{0,\beta}(Y)$ and formula (\ref{thm:tadelg1})   for the tangential gradient of $\int_YK(\cdot,y)\mu(y)\,d\sigma_y$ holds true.
The continuity of the pointwise product in $C^{0,\beta}(Y)$ ensures that the map from 
\[
 C^{0,\beta} (Y)\qquad\text{to}\qquad C^{0,\beta}(Y)\,,
\]
 which takes $\mu$ to $\mu(\cdot) {\mathrm{grad}}_{Y}\int_YK(\cdot,y)\,d\sigma_y$ is linear and continuous. Then   formula (\ref{thm:tadelg1}) implies that 
that map from 
\[
 C^{0,\beta} (Y)\qquad\text{to}\qquad C^{0,\beta}(Y)\,,
\]
 which takes $\mu$ to the function  $\int_Y[{\mathrm{grad}}_{Y,x}K(x,y)](\mu(y)-\mu(x))\,d\sigma_y$ is linear and continuous. Since the kernel $Z(x,y)\equiv  {\mathrm{grad}}_{Y,x}K(x,y)$ satisfy the assumptions of 
  Proposition \ref{prop:mfbd} (i), we conclude that inequality (\ref{eq:maxfubo}) holds true for $X=Y$, 
$Z(x,y)\equiv  {\mathrm{grad}}_{Y,x}K(x,y)$ and thus the proof is complete.\hfill  $\Box$ 

\vspace{\baselineskip}

 \noindent
{\bf Acknowledgement}
The author  acknowledges  the support of the Gruppo Nazionale per l'Analisi Matematica, la Probabilit\`a e le loro Applicazioni   of the Istituto Nazionale di Alta Matematica.\par


\begin{thebibliography}{99}

\bibitem{Cha84}
 I.~Chavel, {\em  Eigenvalues in Riemannian geometry}.
Including a chapter by Burton Randol. With an appendix by Jozef Dodziuk. Pure and Applied Mathematics, 115. Academic Press, Inc., Orlando, FL, 1984.

 

\bibitem{DaLaMu21}
M. Dalla Riva, M.~Lanza de Cristoforis  and P.~Musolino, {\em Singularly Perturbed Boundary Value Problems. A Functional Analytic Approach},  Springer, Cham,  2021.

\bibitem{DoLa17}
F.~Dondi and M.~Lanza de Cristoforis, {\em  Regularizing properties
 of the double layer potential
of  second order elliptic differential operators},   Mem. Differ. Equ.
Math. Phys. 71 (2017), 69--110.

\bibitem{Du76}
J.~Dugundji,  {\em Topology}.
 Allyn and Bacon, Inc., Boston, Mass.-London-Sydney, 1978.
 Reprinting of the 1966 original, Allyn and Bacon Series in Advanced
  Mathematics.

\bibitem{EdKoMe02}
D.E. Edmunds, V. Kokilashvili and A.   Meskhi, {\em Bounded and compact integral operators}, Kluwer, Dordrecht, 2002.

\bibitem{Fo99}
G.B.~Folland, {\em Real analysis. Modern techniques and their applications}. Second edition. John Wiley \& Sons, Inc., New York, 1999.


\bibitem{GaGa04}
J. Garc\'{\i}a-Cuerva and A.E. Gatto,   {\em Boundedness properties of fractional integral operators associated to non-doubling measures}. Studia Math. 162(3) (2004), 245--261.

 \bibitem{GaGa05}
J. Garc\'{\i}a-Cuerva and A.E. Gatto,  {\em Lipschitz spaces and Calderón-Zygmund operators associated to non-doubling measures}. Publ. Mat. 49 (2005), no. 2, 285–296. 	

\bibitem{Gat06}
A.E. Gatto, {\em On fractional calculus associated to doubling and non-doubling measures. Harmonic analysis}, 15–37, Contemp. Math., 411, Amer. Math. Soc., Providence, RI, 2006. 

 \bibitem{Gat09}
A.E.~Gatto, {\em Boundedness on inhomogeneous Lipschitz spaces of fractional integrals singular integrals and hypersingular integrals associated
to non-doubling measures}. Collect. Math. 60, 1 (2009), 101–114.

 


\bibitem{KiHe15}
A.~Kirsch and  F.~Hettlich, {\em  The Mathematical Theory of Time-Harmonic Maxwell's Equations; Expansion-, Integral-, and Variational Methods},  Springer, 2015.  

 
\bibitem{La22a}
 M.~Lanza de Cristoforis, {\em Integral operators in H\"{o}lder spaces  on upper Ahlfors regular sets}, to appear in Atti Accad. Naz. Lincei Rend. Lincei Mat. Appl.  (2023).
 
 \bibitem{La22b}
M.~Lanza de Cristoforis, {\em Classes of kernels and continuity properties of the tangential gradient of an integral operator in H\"{o}lder spaces on a manifold}, submitted, (2022).

 
 
 



\end{thebibliography}
\end{document}